\documentclass[11pt]{article}
\usepackage{amsmath,amsthm,amsfonts,latexsym,amssymb,enumerate,color}
\usepackage{graphicx}

\def\CC{\mathbb C}
\def\DD{\mathbb D}

\def\G{\mathcal G}
\def\GH{\mathcal H}
\def\GK{\mathcal K}
\def\GL{\mathcal L}
\def\GR{\mathcal R}
\def\RR{\mathbb R}
\def\TT{\mathbb T}

\def\diag{\mathop{\rm diag}\nolimits}

\def\ker{\mathop{\rm ker}\nolimits}
\def\dist{\mathop{\rm dist}\nolimits}
\def\ran{\mathop{\rm ran}\nolimits}
\def\limnt{\lim_{z \to \lambda \hbox{ n.t.}}}

\newcommand{\simast}{%
    \ensuremath{%
       \stackrel{\mathsf{\ast}}{\sim}}}

\newtheorem{thm}{Theorem}[section]
\newtheorem{prop}[thm]{Proposition}
\newtheorem{defn}[thm]{Definition}
\newtheorem{lem}[thm]{Lemma}
\newtheorem{cor}[thm]{Corollary}
\newtheorem{rem}[thm]{Remark}

\numberwithin{equation}{section}

\def\beginpf{\begin{proof}}
\def\endpf{\end{proof}}
\def\beq{\begin{equation}}
\def\eeq{\end{equation}}
\def\ds{\displaystyle}
\def\ol{\overline}

\def\Dk{\Delta_{\lambda}} 
\def\abar{\ol{ \ol{A_-}(0) } }

\begin{document}

\title{Dual-band general Toeplitz operators}

\author{M.~Cristina C\^amara,\thanks{
Center for Mathematical Analysis, Geometry and Dynamical Systems,
Instituto Superior T\'ecnico, Universidade de Lisboa, 
Av. Rovisco Pais, 1049-001 Lisboa, Portugal.
 \tt ccamara@math.ist.utl.pt} 
 \and Ryan O'Loughlin\thanks{School of Mathematics, University of Leeds, Leeds LS2~9JT, U.K. \tt mm12rol@leeds.ac.uk}
 \ and  Jonathan R.~Partington\thanks{School of Mathematics, University of Leeds, Leeds LS2~9JT, U.K. {\tt j.r.partington@leeds.ac.uk}}\  
}

\maketitle

\begin{abstract}
We relate dual-band   general Toeplitz operators to 
block truncated Toeplitz operators and, via equivalence after extension, with
Toeplitz operators with $4 \times 4$ matrix symbols.
We discuss their norm, their kernel, Fredhomlness, invertibility and spectral properties
in various situations, focusing on the spectral properties of the
dual-band   shift, which turns out to be considerably
complex, leading to new and nontrivial connections with the boundary behaviour of the associated inner function.
\end{abstract}




\section{Introduction}

Multiband spaces occur naturally in applications.
First, multiband signals are seen in speech processing (see \cite{ABM,BD89,BK,JRH} 
 for example), as an alternative to the Paley--Wiener space
$PW(b)$ of  inverse Fourier transforms of functions in $L^2(-b,b)$, when both high and low frequencies
are to be ignored.
Second, multiplex signal transmission, as a way of sending several signals down the same channel,
has many practical applications, and we refer the reader to \cite{lotsofnames} for a detailed history of the subject
with 137 references, tracing the analysis back to work of Raabe and Shannon in the 1930s and 1940s. Furthermore, in recent years, dual-band   filters have become key components in ubiquitous wireless communication devices such as cellular phones (\cite{JD, LT}).

To see a basic example of a multiband space, which in this case is a dual-band   space, 
choose $0 < a < b$ and consider the inverse Fourier transform of the
space $L^2((-b,-a) \cup (a,b))$, which is a space $M \subset L^2(\RR)$.
Indeed $M=PW(b) \ominus PW(a)$.

If we define the inverse Fourier transform formally by
\[
\hat f(s) = \frac{1}{2\pi}\int_{-\infty}^\infty f(t) e^{ist} \, dt,
\]
then this extends to an isomorphism between $L^2(0,\infty)$ and the Hardy space $H^2(\CC^+)$
on the upper half-plane, and $L^2(0,b-a)$ corresponds to the model space $K_\theta:=H^2(\CC^+) \ominus \theta H^2(\CC^+)$
where $\theta$ is the inner function $\theta(s)=e^{i(b-a)s}$.

It is now clear that the space $M$ above has the orthogonal decomposition 
\beq\label{eq:m}
M=\phi K_\theta \oplus  \psi K_\theta,
\eeq
where 
\beq\label{eq:m2}
\phi(s)=e^{-ibs} \qquad \hbox{and} \qquad \psi(s)=e^{ias}.
\eeq

More generally, let $\theta$ be an inner function in $H^\infty(\CC^+)$, and $\phi$, $\psi$ unimodular
functions in $L^\infty(\RR)$ such that $\phi K_\theta \perp \psi K_\theta$. Then we shall consider  dual-band   spaces $M:= \phi K_\theta \oplus \psi K_\theta$ and general Wiener-Hopf operators (\cite {DS, Speck}) on those spaces. 

Toeplitz operators are a basic example of a general Wiener-Hopf operator where the operator on $L^2$, of the circle or the real line, which is compressed (in this case, to the Hardy space $H^2$) is a multiplication operator; in the latter case we say that we have a general Toeplitz operator. 

Truncated Toeplitz operators of the form $A_g^\theta u=P_{\theta} (gu)$,
where $P_\theta$ denotes the orthogonal projection onto $K_\theta$, are another example of general Toeplitz operators. They
have been much studied, since being formally defined
by Sarason \cite{sarason07}, although they occur much earlier, for example in \cite{AC70,sarason}. Some  recent surveys on the subject are in \cite{CFT16,GR13}.

The dual-band general Toeplitz operator  (abbreviated to dual-band Toeplitz operator in what follows) with symbol $g\in L^\infty$, $T^M_g$, is
defined on the space $M:= \phi K_\theta \oplus \psi K_\theta$ by 
\[
T^M_g u = P_M (gu),
\]
where $P_M$ is the orthogonal projection onto $M$.

In the particular case of the two-interval example given earlier, these are unitarily equivalent to convolution operators restricted
to the union of two intervals. 
Clearly, the same definition can also be made in the more usual situation of $H^2(\DD)$,  except that now $\phi,\psi$ are unimodular in $L^\infty(\TT)$ and $g \in L^2(\TT)$.

Two degenerate cases may appear: the decomposition $M=\overline\theta K_\theta \oplus K_\theta$
gives  the Paley--Wiener space as a special case, and the
decomposition $K_{\theta^2}= K_\theta \oplus \theta K_\theta$ is also a special case. In this paper we will assume that $\bar\phi\psi$ and $\phi\bar\psi$ are not constant multiples of the inner function $\theta$ to avoid those limit cases. 

Note that, in contrast with Toeplitz operators and truncated Toeplitz operators, dual-band Toeplitz operators do not act on spaces of holomorphic functions and $M$ is not a direct sum of model spaces, unless $\phi,\,\psi$ are constant.


In this paper we relate dual-band Toeplitz operators to block truncated Toeplitz operators of a particular form, allowing for the information on the unimodular functions $\phi,\,\psi$ and the symbol $g$ to be encoded in the various components of  the matrix symbol of a block truncated Toeplitz operator. 
By using the
concept of equivalence after extension, these are in their turn related to Toeplitz operators with $4 \times 4$ matrix symbols.
We are thus able to discuss their norm, their kernel, Fredholmness, invertibility and spectral properties
in various important situations. 

We consider in particular the spectral properties of the dual-band
shift. In fact this is a very natural particular case to consider since it has the simplest possible symbol apart from constants and it is an important model to understand the different roles played by $\theta$, on the one hand, and by $\phi,\,\psi$, on the other. There are also links with  Volterra operators, as first noted by Sarason \cite{sarason65}.

This study also highlights the importance of using the equivalence after extension between dual-band  Toeplitz operators
and block Toeplitz operators. Indeed, it allows one to use the powerful Riemann--Hilbert method, used in
a variety of mathematical and physical problems \cite{CCMN,CS,Its}, as well as the theory of Wiener--Hopf factorization
\cite{CDR,MP}, in order 
to study kernels, invertibility and Fredholmness of dual-band Toeplitz operators, to obtain explicit expressions for the inverse operators and the resolvent operators and, for the first time, 
to investigate certain boundary properties of inner functions, such as the existence of an angular derivative
in the sense of Carath\'eodory for $\theta$ in terms of an $L^2$ factorization (\cite{LS}).

Our results are presented in the context of $L^2$ and $L^\infty$ spaces on
the circle $\TT$, or $H^2$ and $H^\infty$ spaces
on the disc $\DD$, but they apply also to the case of $L^2$ and $L^\infty$ spaces on the real line $\RR$
and $H^2$ and $H^\infty$ spaces on the upper half-plane $\CC^+$.

\section{General Toeplitz operators on dual-band   spaces}
\subsection{The dual-band   space $M$}

\begin{prop}\label{prop:21}
Let $\theta \in H^\infty$ be inner, and let $K_\theta= H^2 \ominus \theta H^2$ be the corresponding model
space. Then for $\phi,\psi \in L^\infty $ unimodular, the spaces $\phi K_\theta$ and $\psi K_\theta$ are
orthogonal 
if and only if the truncated Toeplitz operator $A^\theta_{\overline\phi \psi}$ is the zero operator.
\end{prop}

\beginpf
If $f=\psi k_1$ and $g=\phi k_2$ with $k_1, k_2 \in K_\theta$, then
\[
\langle f,g \rangle = \langle \psi k_1,\phi k_2 \rangle  = \langle \overline \phi \psi k_1,   k_2 \rangle 
= \langle P_\theta \overline \phi \psi k_1,   k_2 \rangle 
\]
and this is zero  for all $f,g$ of that form  if and only if $A^\theta_{\overline\phi \psi}=0$.
\endpf

Let now $M:=\phi K_\theta \oplus^\perp \psi K_\theta$, where $\theta$ is inner, $\phi,\psi \in L^\infty $ 
are unimodular and we assume that 
 $A^\theta_{\overline\phi \psi}=0$, i.e., $\overline\phi\psi \in \theta H^2 + \overline\theta \overline{H^2}$
 \cite[Thm. 3.1]{sarason07}. We also assume throughout the paper that $\phi\bar\psi$ and $\bar\phi\psi$ are not constant multiples of $\theta$. 

$M$ is a closed subspace of $L^2$ and  the operator $P_M$ defined by
 \beq\label{eq:n2.1}
 P_M f = \phi P_\theta \overline\phi f + \psi P_\theta \overline \psi f \qquad (f \in L^2)
 \eeq
 is the orthogonal projection from $L^2$ onto $M$.

One  can define a  conjugation $C_M$ on $L^2$ which keeps $M$ invariant (and is therefore a
conjugation on $M$ when restricted to that space). Recall that a conjugation $C$ on a complex Hilbert space $\GH$
is an antilinear isometric involution, i.e.,
\[
C^2= I_\GH \qquad \hbox{and} \qquad \langle Cf,Cg \rangle  = \langle g,f \rangle \quad \hbox{for all} \quad f,g \in \GH.
\]
The study of conjugations, which generalize complex conjugation,
is motivated by applications in physics, in connection with the study
of {\em complex symmetric operators} \cite {GPP, GP1, GP2}. These are the operators $ A \in \GL(\GH)$ such that
$CAC=A^*$ for a conjugation $C$ on $\GH$. We can define a natural
conjugation  $C_\theta$ on any model space $K_\theta$ by 
\beq
C_\theta f = \theta \overline z \overline f,
\eeq
and it is known that any bounded truncated Toeplitz operator is $C_\theta$-symmetric
\cite[Chap. 8]{GMR}.
For $\GH=M$, we have the following result.

\begin{prop}\label{prop:CMconj}
The antilinear operator $C_M$ defined by
\beq
C_M f = \theta \phi \psi \overline z \overline f = \phi \psi C_\theta f
\eeq
is a conjugation on $L^2$ preserving $M$ as an invariant subspace.
\end{prop}
\beginpf
If $f \in M$ has the form $f=\phi k_1+\psi k_2$ with $k_1,k_2 \in K_\theta$, then
\beq 
C_M f =  \psi C_\theta k_1 + \phi C_\theta k_2 \in M,
\eeq
from which the conjugation properties are easily verified.
\endpf

\subsection{A matrix representation }

 Let now $T^M_g $ for  $g \in L^2$ be the operator densely defined in $M$ by
\[
T^M_g f = P_M gf \qquad (f \in L^\infty \cap M),
\]
the density of $L^\infty \cap M$ in $M$ following easily from the density of $L^\infty \cap K_\theta$ in $K_\theta$,
which was given in \cite{sarason07}.

If this operator is bounded, we also denote by $T^M_g$ its unique bounded extension to $M$. The
operator $T^M_g$ is bounded, in particular, whenever $g \in L^\infty$. It is easy to see that
$(T^M_g)^*= T^M_{\overline g}$.

\begin{thm}\label{thm:eqblock}
Let $T^M_g$ be a bounded general Toeplitz  operator on the dual-band   space
$M:=\phi K_\theta \oplus^\perp \psi K_\theta$, where $\theta$ is inner and $\phi,\psi \in L^\infty $ 
are unimodular. Then $T^M_g$ is unitarily equivalent to the block truncated Toeplitz
operator
\beq\label{W}
W=\begin{pmatrix}
A^\theta_g & A^\theta_{\overline\phi \psi g} \\ A^\theta_{\overline\psi \phi g} & A^\theta_g 
\end{pmatrix},
\eeq 
on $K_\theta \oplus K_\theta$. Hence $T^M_g=0$ if and only if each of the 
four truncated Toeplitz operators  composing $W$ is $0$.
\end{thm}

\beginpf
Let $M_\phi$ denote the operator of multiplication by $\phi$,
and similarly for other multiplication operators.
We have the factorization
\begin{eqnarray}\label{eq:3fact}
T^M_g &=&\begin{pmatrix}
M_\phi A^\theta_g M_{\overline \phi} & M_\phi A^\theta_{\overline\phi \psi g}M_{\overline\psi} \\ 
M_\psi A^\theta_{\overline\psi \phi g} M_{\overline \phi} & M_\psi A^\theta_g M_{\overline\psi}
\end{pmatrix} \nonumber \\
&=& \begin{pmatrix}
M_\phi & 0 \\ 0 & M_\psi
\end{pmatrix}
\begin{pmatrix}
A^\theta_g & A^\theta_{\overline\phi \psi g} \\ A^\theta_{\overline\psi \phi g} & A^\theta_g 
\end{pmatrix}
\begin{pmatrix}
M_{\overline \phi} & 0 \\ 0 & M_{\overline\psi}
\end{pmatrix}\,. 
\end{eqnarray}
This has the form $U^*AU$, where $U$ is a unitary operator from $M$ onto $K_\theta \oplus K_\theta$.
Using the fact that $M_{\overline\phi}$ maps $\phi K_\theta$ bijectively to $K_\theta$ and $P_{\phi K_\theta}=M_\phi P_{K_\theta} M_{\overline\phi}$,
it is easy to verify that the identity \eqref{eq:3fact} holds.

\endpf

There are some simplifications possible here, 
since some of these four blocks may be 0. The basic properties of matrix-valued truncated Toeplitz operators were studied in \cite{KT}.





\begin{thm}
If $T^M_g$ is bounded, then it is $C_M$-symmetric.
\end{thm}
\beginpf
We wish to check the identity $T^M_g C_M= C_M (T^M_g)^*$. Note that by 
Proposition \ref{prop:CMconj}, we have 
\[
UC_MU^* = \begin{pmatrix}0 & C_\theta \\ C_\theta & 0 \end{pmatrix},
\]
where $U$ is the unitary mapping given in \eqref{eq:3fact}.
Hence 
\[
U T^M_g   C_M U^* =  \begin{pmatrix}
A^\theta_g & A^\theta_{\overline\phi \psi g} \\ A^\theta_{\overline\psi \phi g} & A^\theta_g 
\end{pmatrix}
\begin{pmatrix}0 & C_\theta \\ C_\theta & 0 \end{pmatrix} = 
 \begin{pmatrix}
A^\theta_{\overline\phi \psi g}  C_\theta & A^\theta_g C_\theta \\ A^\theta_g C_\theta & A^\theta_{\overline\psi \phi g} C_\theta
\end{pmatrix}
\]
while
\[
U C_M ( T^M_g)^*    U^*  =  
\begin{pmatrix}0 & C_\theta \\ C_\theta & 0 \end{pmatrix}
\begin{pmatrix}
A^\theta_{\overline{g}} &  A^\theta_{\psi\overline \phi \overline g}\\ A^\theta_{\phi \overline\psi\overline g}  & A^\theta_
{\overline{g }}
\end{pmatrix} = \begin{pmatrix}
C_\theta A^\theta_{\phi \overline\psi\overline g} & C_\theta A^\theta_{\overline{g}}  \\
C_\theta A^\theta_{\overline{g }}
  &  C_\theta A^\theta_{\psi\overline \phi \overline g} 
\end{pmatrix} ,
\]
and these are equal since $AC_\theta=C_\theta A^*$ for any truncated Toeplitz operator $A$.
\endpf


\section{Equivalence after extension}

\begin{defn} \cite{BTsk, HR, Speck1,dan} The operators $T: X\rightarrow \widetilde{X}$ and $S: Y\rightarrow \widetilde{Y}$ are said to be \emph {(algebraically and topologically) equivalent} if and only if $T=ESF$ where $E,F$ are invertible operators. 
More generally, $T$ and $S$ are \emph{equivalent after extension} if and only if there exist (possibly trivial) Banach spaces $X_0$, $Y_0$, called \emph{extension spaces}, and invertible bounded linear operators $E:\widetilde{Y}\oplus Y_0\rightarrow \widetilde{X}\oplus X_0$ and  $F:X\oplus X_0\rightarrow Y\oplus Y_0$, such that
\beq \label{III.1}
\left(
  \begin{array}{cc}
    T & 0 \\
    0 & I_{X_0} \\
  \end{array}
\right)=E \left(
  \begin{array}{cc}
    S & 0 \\
    0 & I_{Y_0} \\
  \end{array}
\right) F.
\eeq

In this case we say that $T\simast S$.
\end{defn}

It was shown in \cite{CP17} that for $g \in L^\infty$ the scalar Toeplitz operator $A^\theta_g$
is equivalent by extension to the block Toeplitz operator with symbol
\[
\begin{pmatrix}
\overline\theta & 0 \\ g & \theta
\end{pmatrix}.
\]
This result was used in \cite{CP16} to study spectral properties of $A^\theta_g$ and, more generally, to study asymmetric
truncated Toeplitz operators.

%

Motivated by the result of Theorem \ref{thm:eqblock}, we now consider the truncated Toeplitz operator $A^\theta_G$ acting on $K_\theta \oplus K_\theta$, where
$G=\begin{pmatrix}g_{11} & g_{12} \\ g_{21} & g_{22} \end{pmatrix}\in (L^\infty)^{2\times 2}$, and link it with the
Toeplitz operator $T_\G$ acting on $(H^2)^4$, where
\beq\label{eq:bigG32}
\G = \begin{pmatrix} \overline\theta & 0 & 0 & 0 \\ 0 & \overline\theta & 0 & 0 \\
g_{11} & g_{12} & \theta & 0 \\ g_{21} & g_{22} & 0 & \theta
\end{pmatrix}.
\eeq
Clearly,
for $p,q,r,s \in H^2$, we have $(p,q,r,s) \in \ker T_\G$ if and only if
$p,q \in K_\theta$ and 
$\ds \begin{pmatrix}g_{11} & g_{12} \\ g_{21} & g_{22} \end{pmatrix}\begin{pmatrix}p \\ q \end{pmatrix}
+ \theta \begin{pmatrix}r \\ s \end{pmatrix}  \in \overline{H^2_0} \oplus \overline{H^2_0}$.
So $(p,q) \in \ker A^\theta_G$, and likewise given $(p,q) \in \ker A^\theta_G$ there exist $r,s \in H^2$
with $(p,q,r,s) \in \ker T_\G$.\\

%
%


The following theorem shows that the result in \cite[Thm. 2.3]{CP16} can be extended to
block truncated Toeplitz operators and in fact we can give the result more generally for
$n \times n$ blocks. We shall write $P_\theta$ for the orthogonal projection
from $(H^2)^n$ onto $(K_\theta)^n$, and $Q_\theta$ for the complementary projection
from $(H^2)^n$ onto $\theta (H^2)^n$.

\begin{thm}\label{thm:EE13}
Let $G \in (L^\infty)_{n \times n}$ and let $I_n$ be the $n \times n$ identity matrix.
The operator $A^\theta_G=P_\theta G P_\theta :K_\theta^n \to K_\theta^n$ is equivalent after
extension to $T_\G:(H^2)^{2n} \to (H^2)^{2n}$ with 
\beq\label{eq:bigG33}
\G = \begin{pmatrix}
\overline\theta I_n & 0 \\ G & \theta I_n 
\end{pmatrix}.
\eeq
\end{thm}
\beginpf
We have, following the proof of \cite[Thm. 2.3]{CP16}, 
\[
A^\theta_G \simast P_{\theta}G P_{\theta} + Q_{\theta}
\]
because 
\[
\begin{pmatrix}
A^\theta_G & 0 \\ 0 & I_{\theta(H^2)^n}
\end{pmatrix}
= E_1 \begin{pmatrix}
 P_{\theta}G P_{\theta} + Q_{\theta} & 0 \\ 0 & I_{\{0\}^n}
 \end{pmatrix} F_1,
 \]
 where
 
$F_1: K_\theta^n \oplus \theta(H^2)^n \to (H^2)^n \oplus \{0\}^n$ and
$E_1: (H^2)^n \oplus \{0\}^n \to K_\theta^n \oplus \theta(H^2)^n$ are invertible operators, defined in the obvious way.
On the other hand, it is clear that, denoting by $P^+$ the orthogonal projection from $L^2$ onto $H^2$,
\[
P_{\theta}G P_{\theta} + Q_{\theta} \simast 
\begin{pmatrix}
P_{\theta}G P_{\theta} + Q_{\theta} & 0 \\ 0 & P^+ 
\end{pmatrix}: (H^2)^{2n} \to (H^2)^{2n}.
\]
Now,
\[
\begin{pmatrix}
P_{\theta}G P_{\theta} + Q_{\theta} & 0 \\ 0 & P^+
\end{pmatrix}
  = 
\underbrace{\begin{pmatrix}
T_{\theta I_n}- P_\theta G T_{\theta I_n} & P_\theta \\
-P^+ & T_{\overline\theta I_n}
\end{pmatrix}}_{E}
T_\G
\underbrace{\begin{pmatrix}
P^+ & 0 \\
T_{\overline\theta I_n}(P^+ - T_\G) & P^+
\end{pmatrix}}_{F},
\]
where
$E, F: (H^2)^n \to (H^2)^n$ are invertible operators
with 
\[
E^{-1} = \begin{pmatrix} 
T_{\overline \theta I_n} & 0 \\
P^+ + P_\theta G Q_\theta & T_{\theta I_n}
\end{pmatrix}
\]
and
\[
F^{-1} = \begin{pmatrix}
P^+ & 0 \\
-T{\overline\theta I_n}(P^+-T_G) & P^+
\end{pmatrix}.
\]
\endpf

\begin{cor}\label{cor:n2.7}
For $g \in L^\infty$, one has $T^M_g \simast T_\G$ with
\beq\label{eq:slashG}
\G = \begin{pmatrix}
\overline\theta & 0 & 0 & 0 \\
0 & \overline \theta & 0 & 0 \\
g& g\overline\phi \psi & \theta & 0 \\
g\phi\overline\psi & g & 0 & \theta 
\end{pmatrix}.
\eeq
\end{cor}

\beginpf
This is an immediate consequence of Theorems \ref{thm:eqblock} and \ref{thm:EE13}.
\endpf

We clearly have the following corollary of the above.

\begin{cor}
The operators $T^M_g$ and $W$ are invertible (resp., Fredholm) if and only
if $T_\G$ is invertible (resp., Fredholm), with $\G$ given by \eqref{eq:slashG}.
\end{cor}

More general results will be proved later.
%
%

\section{Kernels, ranges and solvability relations}

The equivalence after extension proved in Theorem \ref{thm:EE13} implies certain relations between
the kernels, the ranges, and the invertibility and Fredholm properties of the two operators
$T^M_g$ and $T_\G$ with $\G$ given by \eqref{eq:slashG} (\cite{BTsk, Speck1}), and therefore it implies certain relations between the
solutions of
\beq\label{eq:n3.1}
T_g^M f_M = h_M \qquad \hbox{for a given} \quad h_M \in M
\eeq
and those of
\beq\label{eq:n3.2}
T_\G F_+ = H_+ \qquad \hbox{for a given} \quad H_+ \in (H^2)^4.
\eeq
In this section we study these relations, which also allow for a better
understanding of the equivalence after extension obtained in the previous section.

\begin{thm}\label{thm:n3.1}
$T^M_g f_M=h_M$ with $f_M,h_M \in M$ if and only if
$T_\G F_+=H_+$ with $F_+ = (f_{j+}) \in (H_2^+)^4\,,\,H_+\in (H_2^+)^4$ given by
\begin{eqnarray}\label{eq:n3.8}
f_{1+}=f_{1\theta}=P_\theta \ol\phi f_M, && 
f_{3+}=-P^+\ol\theta g (P_\theta \ol\phi f_M + \ol \phi\psi P_\theta \ol\psi f_M), \nonumber \\
 f_{2+}=f_{2\theta}= P_\theta \ol\psi f_M, &&
f_{4+}=-P^+\ol\theta g (\phi\ol\psi P_\theta \ol\phi f_M +   P_\theta \ol\psi f_M).
\end{eqnarray}
and
\beq\label{eq:n3.7}
 H_+=(0,0,h_{1\theta},h_{2\theta}).
\eeq
Consequently, if $f_M \in \ker T^M_g$, then $F_+ \in \ker T_\G$.
\end{thm}

\begin{proof}
First note that the equation \eqref{eq:n3.2} is equivalent to the Riemann--Hilbert problem
\beq
\G F_+ = F_- + H_+, \qquad F_\pm \in (H_2^\pm)^4,
\eeq
where $H_2^+=H^2$ and $H_2^-= \overline{H^2_0}$.
Let 
\begin{eqnarray}
f_{1\theta}= P_\theta \overline\phi f_M, && h_{1\theta}=P_\theta \overline\phi h_M, \nonumber\\
f_{2\theta}= P_\theta\overline\psi f_M, && h_{2\theta}=P_\theta \overline \psi h_M,
\end{eqnarray}
so that $f_M=\phi f_{1\theta}+\psi f_{2\theta}$ and
$h_M= \phi h_{1\theta}+\psi h_{2\theta}$.
Given $h_M \in M$ we can write, by \eqref{eq:n2.1}, 
$T^M_g f_M=h_M$ if and only if $P_M\,( g\, (\phi f_{1\theta} + \psi f_{2\theta}))= \phi h_{1\theta}+ \psi h_{2\theta}$,
or equivalently,
\[
\phi P_\theta \ol\phi(g\phi f_{1+} + g\psi f_{2+})
+ \psi P_\theta \ol\psi(g\phi f_{1+} + g\psi f_{2+})= \phi h_{1\theta}+\psi h_{2\theta}\,,\,\,\,f_{j+}= f_{j\theta} \in K_\theta.
\]
Since $\phi K_\theta \perp \psi K_\theta$, this
is equivalent to
\begin{eqnarray}
\ol\theta f_{1+}=f_{1-} \in H_2^-, && P_\theta (gf_{1+} + g\ol\phi\psi f_{2+})=h_{1\theta} ,\nonumber \\
\ol \theta f_{2+} = f_{2-} \in H_2^-, && P_\theta (g\phi\ol\psi f_{1+} + g f_{2+})=h_{2\theta},
\end{eqnarray}
which, in its turn, is equivalent to
\[
\bar\theta F_1^+=F_1^-\,\,,\,\,P_\theta \,g\, G_1F_1^+=H_\theta
\]
where 
\beq\label{4.7a}
F_{1\pm}=(f_{1\pm}, f_{2\pm})\,\,,\,\, G_1=\begin{pmatrix}
1 & \bar\phi\psi \\ \phi\bar\psi & 1 
\end{pmatrix}\,\,,\,\,H_\theta=(h_{1\theta}, h_{2\theta}).
\eeq
Equivalently, there exist $F_{2\pm}\in H_{2\pm}^2$ such that 
\beq\label{4.5a}
\bar\theta\,F_{1+}=F_{1-}\,\,\,,\,\,\,g\,G_1F_{1+}=H_\theta+F_{2-}-\theta F_{2+}.
\eeq
This system determines $F_{2+}$ in terms of $F_{1+}$ as $F_{2+}=P^+(\bar\theta g\,G_1F_{1+})$. So, identifying $F_\pm\in (H^2_\pm)^4$ with $(F_{1+},F_{2+})$, \eqref{4.5a} is equivalent to
\[\G F_+ = F_- + H_+\]
with $H_+=(0,0,h_{1\theta},h_{2\theta})$ and $F_+=(f_{j+})$ where
\begin{eqnarray}
f_{1+}=P_\theta \ol\phi f_M, && 
f_{3+}=-P^+\ol\theta g (P_\theta \ol\phi f_M + \ol \phi\psi P_\theta \ol\psi f_M), \nonumber \\
 f_{2+}= P_\theta \ol\psi f_M, &&
f_{4+}=-P^+\ol\theta g (\phi\ol\psi P_\theta \ol\phi f_M +   P_\theta \ol\psi f_M).
\end{eqnarray}

\end{proof}


\begin{thm}\label{thm:n3.2}
$T_\G F_+=H_+$ with $F_+=(\tilde f_{j+}),H_+ =(h_{j+})\in (H^2)^4$,
if and only if
$T^M_g f_M = h_M$, where
\beq\label{eq:nfM}
f_M=\phi(\tilde f_{1+}-\theta h_{1+})+\psi(\tilde f_{2+}-\theta h_{2+})
\eeq
and
\beq\label{eq:nhM}
h_M=\phi P_\theta (h_{3_+}-g\theta (h_{1+}+\ol\phi h_{2_+}))+
\psi P_\theta(h_{4+}-g\theta (\phi\ol\psi h_{1+} + h_{2+})).
\eeq
\end{thm}

\beginpf
We have $T_\G F_+=H_+$ if and only if $\G F_+=F_- + H_+$, with $F_- \in (H_2^-)^4$.
Let $\tilde F_{1\pm}=(\tilde f_{1\pm},\tilde f_{2\pm})$ and $\tilde F_{2\pm}=(\tilde f_{3\pm},\tilde f_{4\pm})$ where $(\tilde f_{j\pm})=F_\pm$ and let $H_{1+}= (h_{1+},h_{2+})\,,\,H_{2+}=(h_{3+}\,,\,h_{4+}$. Then $\G F_+=F_- + H_+$ if and only if 
\beq
\left\{
\begin{array}{rcl}
\ol\theta \tilde F_{1+} &=& \tilde F_{1-}+ H_{1+}   \\ 
g\,G_1\tilde F_{1+} +\theta \tilde F_{2+} &=& \tilde F_{2-}+H_{2+}
\end{array}
\right.\label{eq:n4.11}
\eeq
where $G_1$ is given in \eqref{4.7a}.

This  in turn is equivalent to the system of equations

\beq
\left\{
\begin{array}{rcl}
\ol\theta (\tilde F_{1+}-\theta H_{1+}) = \tilde F_{1-}   \\ 
g\,G_1(\tilde F_{1+}-\theta H_{1+}) +\theta (\tilde F_{2+}+ P^+ g\,G_1H_{1+}-P^+\bar\theta H_{2+})\\
-P_\theta H_{2+}-P_\theta \theta g\,G_1H_{1+}
= \tilde F_{2-}+ P^-\theta g\,G_1H_{1+}
\end{array}
\right.
\eeq

Taking $ F_{1+}=\tilde F_{1+}-\theta H_{1+}\,,\,F_{2+}=\tilde F_{2+}+ P^+ g\,G_1H_{1+}-P^+\bar\theta H_{2+}$ we get

\beq
\left\{
\begin{array}{rcl}
\ol\theta  F_{1+}= \tilde F_{1-}   \\ 
g\,G_1 F_{1+} +\theta  F_{2+}
-P_\theta H_{2+}-P_\theta \theta g\,G_1H_{1+}\\
= \tilde F_{2-}+ P^-\theta g\,G_1H_{1+}+P_\theta (H_{2+}+ \theta g\,G_1H_{1+})
\end{array}
\right.
\eeq

By Theorem \ref{thm:n3.1} this is equivalent to $T^M_g f_M = h_M$ with $f_M$ and $h_M$ given by
\eqref{eq:nfM}
and
\eqref{eq:nhM}.
\endpf

\begin{cor}\label{4.3}
If $F_+ \in \ker T_\G$, with $F_+=(\tilde f_{j+})\in (H^2)^4$, then $f_M \in \ker T^M_g$ with 
$f_M= \phi \tilde f_{1+} + \psi \tilde f_{2+}$.
\end{cor}

Note that  \eqref{eq:n4.11} shows that any element of the kernel of $T_\G$ is determined
by its first two components $\tilde f_{1+}$ and $\tilde f_{2+}$, since
\[
\tilde f_{3+} = P^+\ol\theta g(\tilde f_{1+} + \ol\phi\psi \tilde f_{2+})
\]
and
\[
\tilde f_{4+} = P^+\ol\theta g(\phi\ol\psi\tilde f_{1+} +   \tilde f_{2+}).
\]
Let $P_{1,2}$ be the projection defined by $P_{1,2}(x,y,u,v)=(x,y)$.
\begin{cor}\label{4.4}
The map \[\GK: \ker T_g^M \to \ker T_\G\,\,,\,\,\,\GK f_M = (f_{1+},f_{2+},f_{3+},f_{4+})\]
with
$f_{j+}$ given by \eqref{eq:n3.8} is an isomorphism.  We have 
\[
\begin{array}{l}
\ker T^M_g = \GK^{-1} \ker T_\G 
= \{\phi f_{1+}+\psi f_{2+}: (f_{1+},f_{2+})\in P_{1,2} \ker T_\G\}.
\end{array}
\]
\end{cor}
From Theorems \ref{thm:n3.1} and \ref{thm:n3.2} we also obtain the following regarding ranges.

\begin{cor}\label{4.5}
With the same notation as above, \\
(i) $(h_{1+},h_{2+},h_{3+},h_{4+}) \in \ran T_\G$ if and only if
\[
\phi P_\theta (h_{3+}-g\theta (h_{1+}+ \ol\phi \psi h_{2+}))
+ \psi P_\theta (h_{4+}-g\theta (\phi\ol\psi h_{1+}+h_{2+}))
\in \ran T^M_g;
\]
(ii) $\phi h_{1\theta}+\psi h_{2\theta} \in \ran T^M_g$ if and only if $(0,0,h_{1\theta},h_{2\theta}) \in \ran T_\G$.
\end{cor}
Moreover, we obtain a relation between the inverses of $T^M_g$ and $T_\G$ when these operators are invertible.

\begin{cor}\label{4.6}
$T^M_g$ is invertible if and only if $T_\G$ is invertible and, in that case,
$(T^M_g)^{-1} = [\phi P_1,\psi P_2, 0, 0] T_\G^{-1} U_0$,
where $P_j(x_1,x_2,x_3,x_4)=x_j$ and
$U_0:M \to (H_2^+)^4$, is given by
$U_0 h_M = (0,0, P_\theta \ol\phi h_M, P_\theta \ol \psi h_M)$.
\end{cor}

We can also relate the kernels of $T_g^M$ and its adjoint as follows.

\begin{thm}\label{4.7}
If $g \in L^\infty$, then $\ker T^M_g \simeq \ker (T_g^M)^* = \ker T^M_{\overline g}$.
\end{thm}

\beginpf
Since $T^M_g \simast T_\G$, we have that
$\ker T^M_g \simeq \ker T_\G$ and $\ker (T^M_g)^* \simeq \ker T_\G^*$. So
it is enough to prove that $\ker T_\G$, with $\G$ given by
\eqref{eq:slashG}, is isomorphic to $\ker T_\G^* = \ker T_{\overline \G^T}$. Since
\[
\ker T_\G= \{ \phi_+ \in (H^2)^4: \G\phi_+ = \phi_- \in (\overline{H^2_0})^4 \},
\]
we have 
\[
\G \phi_+ = \phi_- \iff \overline z \overline{\phi_+} = \overline{\G^{-1}} (\overline z \overline{\phi_-})
\iff \overline{\G}^{-1} \psi_+ = \psi_-,
\]
where $\psi_+= \overline z \overline{\phi_-} \in (H^2)^4$ and 
$\psi_- = \overline z \overline{\phi_+} \in (\overline{H^2_0})^4$.
Since
\begin{eqnarray*}
\overline \G^{-1} &=&
\begin{pmatrix}
\ol\theta & 0 & 0 & 0 \\
0 & \ol \theta & 0 & 0 \\
-\ol g & -\ol g \phi\ol\psi & \theta & 0 \\
-\ol g \ol \phi \psi & - \ol g & 0 & \theta
\end{pmatrix}\\
&=&
\begin{pmatrix}
0 & 0 & 0 & 1 \\
0 & 0 & 1 & 0 \\
0  & -1 & 0 & 0 \\
-1 & 0 & 0 & 0
\end{pmatrix}
\underbrace{\begin{pmatrix}
\theta & 0 & \ol g & \ol g \ol \phi \psi \\
0 & \theta & \ol g \phi \ol\psi & \ol g \\
0 & 0 & \ol\theta & 0 \\
0 & 0 & 0 & \ol\theta
\end{pmatrix}}_{\overline \G^T}
\begin{pmatrix}
0 & 0 & 0 & -1 \\
0 & 0 & -1 & 0 \\
0 & 1 & 0 & 0 \\
1 & 0 & 0 & 0
\end{pmatrix},
\end{eqnarray*}
it is clear that $\ker T_{\overline\G^{-1}} \simeq \ker T_{\overline\G^T}$.
\endpf

\begin{cor}\label{cor:n3.8}
If $T^M_g$ is Fredholm, then it has index $0$. Consequently, 
$T^M_g$ is invertible if and only if it is Fredholm and injective.
\end{cor}

\section{The norm and the spectrum for analytic symbols} 

Clearly the norm of the dual-band  Toeplitz operator $T^M_g$ is the same as the norm of the block truncated Toeplitz operator $W$. 
One important case that can be analysed is when $g$ is
in $H^\infty$ (in the language of dual-band  signals, this corresponds to a causal convolution
on $L^2((-b,-a) \cup (a,b))$).

The following is an easy generalization of scalar results which apparently go back to \cite{sarason}.

\begin{prop}
Suppose that the symbol
\[
\Phi:=\begin{pmatrix}
g & \overline\phi \psi g \\ \overline\psi \phi g & g
\end{pmatrix}
\]
is in $(H^\infty)_{2 \times 2}$. Then 
\[
\|T^M_g\|=\|W\|=
\dist(\overline\theta \Phi, (H^\infty)_{2 \times 2} = 
\|\Gamma_{\overline\theta \Phi}\|,
\]
 where the vectorial Hankel operator $\Gamma_{\overline\theta \Phi} : (H^2)^2 \to (L^2 \ominus H^2)^2$
is defined by $\Gamma_{\overline\theta \Phi}v = P_{(L^2 \ominus H^2)^2} \overline \theta \Phi v$.
\end{prop}

\beginpf
Since the symbol $\Phi$ is analytic, if we write 
 $(u_1,u_2) \in H^2 \oplus H^2$
as $(k_1+\theta \ell_1, k_2 + \theta \ell_2)$ 
with the $k_j $ in $K_\theta$ and $\ell_j \in H^2$, then
we have $W(u_1,u_2)=W(k_1,k_2)$ implying that the norm of the truncated Toeplitz operator $W$ is
the same when  the domain is $H^2 \oplus H^2$ or $K_\theta \oplus K_\theta$. 

But $Wu  = \theta (P_{-} \oplus P_{-}) \overline\theta \Phi u = \theta \Gamma_{\overline\theta \Phi}u$, and so
\[
\|W\|= \|\Gamma_{\overline\theta \Phi}\|= \dist(\overline\theta \Phi, H^\infty(M_2(\CC)),
\]
 by 
the vectorial form of Nehari's theorem~\cite[Sec.\ 2.2]{peller}.\\

\endpf

Some results on the spectrum of $W$  can be derived using known results on the scalar case,
particularly in the context of Proposition~\ref{prop:21}. Note that the hypotheses
of this theorem are satisfied in the original example given by
\eqref{eq:m} and \eqref{eq:m2}.

\begin{thm}
\label{thm:spec}
Suppose that $g \in H^\infty$ and that $\overline\phi \psi \in \theta H^\infty$ or $\overline\phi \psi \in \overline{\theta H^\infty}$. Then
\[
\sigma(T^M_g)= \sigma(A^\theta_g)= 
\{ \lambda \in \CC: \inf_z( |\theta(z)|+| g(z)-\lambda|)=0 \},
\]
where the infimum is taken over $\DD$ 
.
\end{thm}

 \beginpf
Assume that $g \in H^\infty$ and that $\overline\phi \psi \in \theta H^\infty$. Note  that $W$ has the form 
\[
W= \begin{pmatrix}
A^\theta_g & 0 \\ A^\theta_{\overline\psi \phi g} & A^\theta_g
\end{pmatrix}
\]
and we claim that $W$ is invertible if and only if $A^\theta_g$ is. 
For the necessity note that for arbitrary block operator matrices, if we have
\[
\begin{pmatrix} A & 0 \\ B & A \end{pmatrix}
\begin{pmatrix} P & Q \\ R & S \end{pmatrix}=
\begin{pmatrix} P & Q \\ R & S \end{pmatrix}
\begin{pmatrix} A & 0 \\ B & A \end{pmatrix}=
\begin{pmatrix} I & 0 \\ 0 & I \end{pmatrix},
\]
then $AP=I$ and $SA=I$, so $S=SAP=P$, and $A$ is invertible with inverse $P$.

The sufficiency follows from the formula
\[
W^{-1} = 
\begin{pmatrix}
(A^\theta_g){}^{-1} & 0 \\ -(A^\theta_g){}^{-1}A^\theta_{\overline\psi \phi g} (A^\theta_g){}^{-1} & (A^\theta_g){}^{-1} 
\end{pmatrix}.
\]
The spectrum of $A^\theta_g$ for $g \in H^\infty$ is described in \cite[p. 66]{nik},
and the $H^\infty(\CC^+)$ case may be found in \cite{CP16}.
\endpf

For the essential spectrum of $A^\theta_g$ we may similarly prove the following result.

\begin{thm}
Suppose that $g \in H^\infty$ and that $\overline\phi \psi \in \theta H^\infty$. Then
\[
\sigma_e(T^M_g)= \sigma_e(A^\theta_g)= \{\lambda \in \CC: \liminf_{z \to \xi}
(|\theta(z)|+|g(z)-\lambda|)=0 \hbox{ for some } \xi \in \TT\},
\]
where $z$ is taken in $\DD$. 

\end{thm}
\beginpf
The method of proof of Theorem \ref{thm:spec} adapted to
inversion modulo the compact operators (i.e., in the Calkin algebra)
shows directly that $\sigma_e(T^M_g)=\sigma_e(A^\theta_g)$,
and an expression for this is known from results in \cite{bessonov,CP16}.
\endpf

%

These results are of particular interest in the case of the {\em restricted shift\/} or
{\em truncated shift\/} $S_M$ on $M$, 
with $g(z)=z$. We thus have

\begin{cor}
\label{corx:5.4}
 If $\ol\phi\psi \in \theta H^\infty$, then for the restricted shift $S_M$ on $M$ we have 
\[
\sigma(S_M)= \{\lambda \in \CC: \inf_{z \in \DD}( |\theta(z)|+| z-\lambda|)=0 \},
\]
and
\begin{eqnarray*}
\sigma_e(S_M)&=&
\{\lambda \in \CC: \liminf_{z \to \xi}
(|\theta(z)|+|z-\lambda|)=0 \hbox{ for some } \xi \in \TT\} \\
&=& \{\lambda \in \TT: \liminf_{z \to \lambda}
(|\theta(z)| )=0\}.
\end{eqnarray*}
\end{cor}






\section{The double-band shift: spectral properties}
\label{sec:6}

In order to have matrix symbols which are essentially bounded, and since $\bar\psi\phi\in \bar\theta\overline{H^2}+\theta H^2$, we assume here that $\ol\psi\phi = A_-\ol\theta + A_+\theta$ with $A_+\in H^\infty$ and $A_- \in \ol{H^\infty}$. In this case, for $g=z-\lambda$, we have  in $W$ (see \eqref {W}) $A^\theta_{\bar\phi\psi g}=A^\theta_{\bar A_+\bar\theta (z-\lambda)}$ and $A^\theta_{\bar\psi\phi g}=A^\theta_{ A_-\bar\theta (z-\lambda)}$.
Using the result of Corollary \ref{cor:n2.7} we thus associate to the operator $T^M_{z-\lambda}$,
with $\lambda\in\CC$, the matrix symbol
\beq\label{6.1}
\G_\lambda = \begin{pmatrix}
\ol\theta & 0 & 0 & 0 \\
0 & \ol\theta & 0 & 0 \\
z-\lambda & (z-\lambda)\overline{A_+}\ol\theta & \theta & 0 \\
(z-\lambda) A_-\ol\theta & z-\lambda & 0 & \theta 
\end{pmatrix}.
\eeq


\subsection{Eigenvalues and eigenspaces}\label{sec6.1}
By Corollaries \ref{4.3} and \ref{4.4} we have that $\ker T_{z - \lambda}^M \simeq \ker T_{\G_\lambda}$ and $f^M \in \ker T_{z- \lambda}^M $ if and only if
\beq\label{6.2}
f^M = \phi f_{1+} + \psi f_{2+}
\eeq
where $f_{1+}$ and $f_{2+}$ are the two first components of $\ker T_{\G_\lambda} \subset (H^2)^4$. Note that $\ker T_{\G_\lambda}$ consists of the solutions $f_{+} \in (H^2)^4$ of 
\beq\label{6.3}
\G_\lambda f_+=f_- ,\,\,\,\ \text{with } f_+ \in (H^2)^4, f_- \in (\overline{H^2_0})^4 .
\eeq
It is easy to see from \eqref{6.3} that all components of $f_+$ and $f_-$ are determined by $f_{1+}$ and $f_{2+}$.

We will consider the cases $\lambda \in \mathbb{D}$, $\lambda \in \mathbb{D}^- = \{ z \in \mathbb{C} : |z| > 1 \} $ and $ \lambda \in \mathbb{T}$ separately.

\begin{prop}\label{prop6.1}
If $ \lambda \in \mathbb{D}$ then:
\newline
(i) $\ker T_{z- \lambda}^M = \{ 0 \} \iff \Dk:= \theta(\lambda)^2 - \ol{A_+(0)}\,\ol{\ol{A_-}(0)}(1-\ol{\theta(0)}\theta(\lambda))^2 \ne 0$
\\

\noindent (ii) $\dim \ker T_{z- \lambda}^M = 1 \iff  \Dk =0$ and $|\theta(\lambda)| + |A_+ (0) | + |\ol{A_-}(0)| \ne 0$
\\

\noindent (iii) $\dim \ker T_{z- \lambda}^M = 2 \iff  \Dk =0$ and $\theta(\lambda) = A_+ (0) = \overline{A_-}(0) =0$
\end{prop}

\beginpf
From \eqref{6.1} and \eqref{6.3} we get, for $f_{\pm} = (
f_{1 \pm}, f_{2 \pm}, f_{3 \pm}, f_{4 \pm} ) \ne 0$,
\beq\label{6.4}
( z- \lambda) \begin{pmatrix}
f_{1+} \\
f_{2+}
\end{pmatrix} 
+ \theta \begin{pmatrix}
f_{3+} \\
f_{4+}
\end{pmatrix} = - (z - \lambda) \begin{pmatrix}
 \overline{A_+} f_{2-} \\
A_- f_{1-} 
\end{pmatrix} +
\begin{pmatrix}
f_{3-} \\
f_{4-}
\end{pmatrix} = \begin{pmatrix}
k_1 \\
k_2
\end{pmatrix}
\eeq
with $k_1, k_2 \in \mathbb{C}$, since the left-hand side is in $(H^2)^2$ and the right-hand side of the first equality is in $\overline{(H^2)}^2$. Therefore
\[
(z-\lambda) f_{1+} + \theta f_{3+} = k_1
\]
and, since $\overline{\theta} f_{1+} = f_{1-}$, we have
\beq\label{6.5}
(z- \lambda) f_{1-} - \overline{\theta} k_1 = - f_{3+} = C_1 \in \mathbb{C},
\eeq
so $(z-\lambda) f_{1+} = k_1 + C_1 \theta$ and we have $k_1 = - C_1 \theta (\lambda)$. It follows that 
\beq\label{6.6}
f_{1+} = C_1 \dfrac{\theta-\theta(\lambda)}{z-\lambda}.
\eeq
Analogously, 
\beq\label{6.7}
f_{2+} = C_2 \dfrac{\theta-\theta(\lambda)}{z-\lambda}
\eeq
and it follows that 
\begin{align*}
f_{1-}&= C_1 \dfrac{1-\ol\theta \theta(\lambda)}{z-\lambda}\,\,,   && f_{2-}= C_2 \dfrac{1-\ol\theta \theta(\lambda)}{z-\lambda} \\
f_{3+}&=-C_1 \,\,, && f_{4+}=-C_2 \\
f_{3-} &= -\theta(\lambda)C_1 + \ol{A_+}(1-\ol\theta \theta(\lambda))C_2\,\,,\\
f_{4-} &= -\theta(\lambda)C_2 + A_-(1-\ol\theta \theta(\lambda))C_1\,\,.
\end{align*}
Since $f_{3-}, f_{4-} \in \overline{H^2_0}$, we must have
\begin{eqnarray}\label{6.8}
 -\theta(\lambda)C_1 + \ol{A_+(0)} (1-\ol{\theta(0)} \theta(\lambda))C_2 &=& 0 \nonumber \\
 -\theta(\lambda)C_2 + \ol{ \ol{A_-}(0) } (1-{\ol\theta(0)} \theta(\lambda))C_1 &=& 0.
 \end{eqnarray}
If $C_1 = C_2 = 0$, we get from \eqref{6.4} that $f_{+} = f_- = 0$. A necessary and sufficient condition for \eqref{6.8} to have non zero solutions, $C_1, C_2$, is that the determinant of the system is zero, i.e., $\Dk = 0$. So $(i)$ holds.

If $\Dk =0$, then \eqref{6.8} is equivalent to 
\begin{equation}\label{6.9}
\theta(\lambda) C_1 = \overline{A_+(0)}  (1-\overline{\theta(0)}\theta(\lambda))  C_2.
\end{equation}
If $ \theta( \lambda) \ne 0$, we must have also $A_+ (0), \overline{A_-}(0) \ne 0$ and the system \eqref{6.8} is equivalent to $ C_1= \ol{A_+(0)}\, \frac{1-\ol{\theta(0)}\theta(\lambda)}{\theta(\lambda)}\, C_2$. If $\theta( \lambda) = 0$ then $A_+ (0) \overline{A_-}(0) = 0$.

If $A_+(0) = 0, \overline{A_-}(0) \ne 0$, then
\[
C_1 = 0,\text{ } f_{1+} = 0 ,\text{ }\,\, f_{2+} = \beta_2 \frac{\theta}{z - \lambda}, \text{ with } \beta_2 \in \mathbb{C};
\]
if $A_+(0) \ne 0, \overline{A_-}(0) = 0$, then
\[
C_2=0,\text{ } f_{1+} = \beta_1 \frac{\theta}{z - \lambda}, \text{ with } \beta_1 \in \mathbb{C}, \,\,f_{2+} = 0;
\]
if $A_+(0) = \overline{A_-}(0)$, we have
\[
f_{1+} = \beta_1 \frac{\theta}{z - \lambda}, \text{ }f_{2+}= \beta_2 \frac{\theta}{z - \lambda} \text{ with } \beta_1, \beta_2 \in \mathbb{C}.
\]
So $(ii)$ and $(iii)$ hold.
\endpf
\begin{prop}\label{prop6.2}
If $\lambda \in \mathbb{D}^-$ then :
\newline
(i) $\ker T_{z- \lambda}^M = \{ 0 \} \iff \widetilde{\Dk} := 1- \overline{A_+ (0)} \overline{\overline{A_-}(0)} ( \overline{\theta(0)} - \overline{\theta(\lambda)})^2 \neq 0$,
\\

\noindent (ii) $\dim \ker T_{z- \lambda}^M = 1 \iff  \widetilde{\Dk} = 0$.
\end{prop}

\beginpf
From \eqref{6.4} we get 
\[
(z - \lambda) f_{1-} - \overline{\theta} k_1 = C_1 \iff (z - \lambda) f_{1-} = C_1 + \overline{\theta} k_1 .
\]
Replacing $ z $ $(\in \mathbb{D}^-)$ by $ \lambda$ we see that
\[
C_1 =  - \overline{\theta}( \lambda) k_1, \text{ with } \overline{\theta} (\lambda) = \overline{\theta(1/\overline{\lambda})}.
\]
Analogously, we get $C_2 = - \overline{\theta}(\lambda) k_2$, so 
\[
f_{1-} = k_1 \dfrac{\overline{\theta}-\overline{\theta}(\lambda)}{z-\lambda},\text{ } f_{2-} = k_2 \dfrac{\overline{\theta}-\overline{\theta}(\lambda)}{z-\lambda}.
\]
It follows that
\begin{align}\label{6.10}
f_{1+}&= k_1 \dfrac{1-\theta \ol\theta(\lambda)}{z-\lambda}\,\,,   && f_{2+}= k_2 \dfrac{1-\theta \ol\theta(\lambda)}{z-\lambda} \nonumber \\
f_{3+}&= \ol\theta (\lambda) k_1 \,\,, && f_{4+}=\ol\theta (\lambda) k_2, \nonumber\\
f_{3-} &= k_1 + \overline{A_+} k_2 ( \overline{\theta} - \overline{\theta}(\lambda))\,\,, && f_{4-} = k_2 + A_- k_1 ( \overline{\theta} - \overline{\theta}(\lambda))\,\,.
\end{align}
Since $f_{3-}, f_{4-} \in \overline{H^2_0}$, we must have, from \eqref{6.10}, 
\begin{eqnarray}\label{6.11}
k_1 + \ol{A_+(0)} k_2 ( \overline{ \theta (0)} - \overline{\theta}(\lambda)) &=& 0  \nonumber \\
k_2 + \overline{\ol{A_-}(0)} k_1 ( \overline{ \theta (0)} - \overline{\theta}(\lambda)) &=& 0.
\end{eqnarray}
and discussing this system as in the proof of Proposition \ref{prop6.1} we conclude that $(i)$ and $(ii)$ hold. Moreover, if $\widetilde{\Dk} = 0$ we see that
\begin{align*}
f_{1+}&= - \beta \overline{A_+ (0)} ( \overline{ \theta (0)} - \overline{\theta}(\lambda)) \dfrac{1-\theta \ol\theta(\lambda)}{z-\lambda}, \\
f_{2+}&= \beta \dfrac{1-\theta \ol\theta(\lambda)}{z-\lambda},
\end{align*}
with $\beta \in \mathbb{C}$, which determines the kernel of $T_{z - \lambda}^M$ by \eqref{6.2}.
\endpf

\noindent To study $\ker T_{z - \lambda}^M$ where $ \lambda \in \mathbb{T}$ we will need the following two lemmas.
\begin{lem}\label{lem6.3}
If $\phi_+ \in H^2$ then $\phi_+(z) (z-\lambda) \to 0$ when
$z \to \lambda \in \TT$ nontangentially in $\mathbb D$.
\end{lem}

\beginpf
For $\phi_+ \in H^2$ we have $|\phi_+(w) | \le \|\phi_+\| \, \|k_w\|$,
where $k_w$ is the reproducing kernel function $k_w(z)=1/(1-\ol w z)$, with
$\|k_w\| = 1/\sqrt{1-|w|^2}$. Hence 
\[
|\phi(z)(z-\lambda)| \le \frac{ |z-\lambda| \, \|\phi_+\|}{\sqrt{1-|z|^2}},
\]
and this tends to $0$ if $z$ tends nontangentially to $\lambda \in \TT$,
as this means that $|z-\lambda| \le C(1-|z|)$ for some constant $C$.
\endpf

We say that an inner function $\theta$ has an angular derivative in the sense of Carathéodory (ADC) if and only if $\theta$ has a nontangential limit
\begin{equation*}
    \theta (\lambda) = \limnt \theta(z)
\end{equation*}
with $| \theta (\lambda)| = 1$ and the difference quotient $\ds \frac{\theta(z)-\theta(\lambda)}{z-\lambda}$ has a nontangential limit $\theta^{'}(\lambda)$ at $\lambda$ ( \cite{GMR, sarason07}). By Theorem 7.4.1 in \cite{GMR}, $\theta$ has an ADC at $\lambda \in \mathbb{T}$ if and only if there exists $a \in \mathbb{T}$ such that $\ds \frac{\theta(z)-a}{z-\lambda} \in H^2$, which implies, by Lemma \ref{lem6.3}, that there exists the limit $\theta( \lambda)$ and we have $\theta(\lambda) = a$. Thus we have:
\begin{lem}\label{lem6.4}
$\theta$ has an ADC at $\lambda \in \mathbb{T}$ if and only if:\\
(i) $\limnt \theta(z)$ exists in $\mathbb{C}$ (denoted by $\theta (\lambda)$) \\
and \\
(ii) $\ds \frac{\theta(z)-\theta(\lambda)}{z-\lambda} \in H^2$.
\end{lem}
\noindent We denote by $\mathbb{T}_{\text{ADC}}$ the set of all $ \lambda \in \mathbb{T}$ where $\theta$ has an ADC. 

\begin{prop}\label{prop6.5}
Let $\lambda \in \mathbb{T}$. If $\ker T_{z- \lambda}^M \ne \{ 0 \}$ then $ \lambda \in \mathbb{T}_{\text{ADC}}$.
\end{prop}
\beginpf
Consider again \eqref{6.3}, now with $ \lambda \in \mathbb{T}$, and assume that $f_+, f_- \ne 0$. then, as in the proof of Proposition \ref{prop6.1}, we get \eqref{6.4} and, if $k_1 = k_2 = 0$, we have that $ ( z - \lambda) f_{1+} + \theta f_{3+} = 0$, which implies that 
\[
( z - \lambda) f_{1-} = - f_{3+} = C \in \mathbb{C}.
\]
It follows that $ f_{1-} = \frac{C}{z- \lambda}$, so $C=0$ and therefore $f_{1-} = f_{1+} = 0$. Analogously, from \eqref{6.4} we get $f_{2+}=0$ if $k_2 = 0$. Thus, to have a non-zero solution to \eqref{6.3} either $k_1$ or $k_2$ must be different from 0. Assume that $k_1 \ne 0$. Then, from \eqref{6.4}, 
\[
(z - \lambda) f_{1+} +  \theta f_{3+}= k_1 \implies (z - \lambda) f_{1-} - \overline{\theta} k_1 = - f_{3+} = C \in \mathbb{C} \implies 
\]
\[
(z - \lambda) f_{1-} = C + \overline{\theta} k_1 \implies (z - \lambda) f_{1+} = C \theta + k_1.
\]
From Lemma \ref{lem6.3} it follows now that there exists $ \theta ( \lambda)$ and $ f_{1+} = C \dfrac{\theta-\theta(\lambda)}{z-\lambda} \in H^2$, so that by Lemma \ref{lem6.4} $ \theta$ has an ADC at $\lambda$.
\endpf
\begin{cor}\label{cor6.6}
If $ \lambda \in \mathbb{T} \setminus \mathbb{T}_{\text{ADC}}$ then $ \ker T_{z - \lambda}^M = \{ 0 \}$.
\end{cor}
\begin{prop}\label{prop6.7}
Let $ \lambda \in \mathbb{T}_{\text{ADC}}$. Then, for $\Dk$ defined as in Proposition \ref{prop6.1}, 
\newline
(i) if $\Dk \ne 0$ then $ \ker T_{z - \lambda}^M = \{ 0 \}$,
\\
\noindent (ii)  if $\Dk = 0$ then $\dim \ker T_{z - \lambda}^M = 1$.
\end{prop}
\beginpf
Analogous to the proof of Proposition \ref{prop6.1}, taking into account Proposition \ref{prop6.5} and noting that now we cannot have $ \theta ( \lambda )= 0$.
\endpf
\noindent We summarise the previous results as follows.
\begin{thm}\label{thm6.8}
(i) $\ker T_{z - \lambda}^M$ and $\ker (T_{z - \lambda}^M)^*$ have the same finite dimension for all $ \lambda \in \mathbb{C}$.
\newline
(ii) $ \lambda \in \sigma_p ( T_z^M)$ if and only if $ \lambda \in \mathbb{D}$, $\Dk = 0$, or $ \lambda \in \mathbb{D}^-$, $ \widetilde{\Dk}=0,$ or $ \lambda \in \mathbb{T}_{\text{ADC}}$, $\Dk = 0 ,$ where $ \Dk $ and $ \widetilde{ \Dk}$ are defined in Proposition \ref{prop6.1} and \ref{prop6.2}, respectively.  
\end{thm}
\begin{rem}
For $\lambda \in \sigma_p ( T_z^M)$, the previous results provide a description of the eigenspace in each case.
\end{rem}

\noindent The following corollary applies in particular to the case considered in Theorem \ref{thm:spec}.
\begin{cor}\label{cor6.10}
If $A_+ (0) \overline{A_-}(0) = 0$ then $ \lambda \in \sigma_p ( T_z^M) \iff \lambda \in \mathbb{D}$, $\theta (\lambda) = 0$.
\end{cor}
\begin{cor}\label{cor6.11}
If $A_+ (0) \overline{A_-}(0) \ne 0$ then
\begin{align}
\lambda& \in \sigma_p(T^M_z) \cap \mathbb{D} \implies \theta(\lambda) = 0\,\,  \nonumber \\
\lambda& \in \sigma_p(T^M_z) \cap \mathbb{D}^- \implies \theta(1/ \overline{\lambda}) \neq 0 \,\, \nonumber
\end{align}
\end{cor}

\subsection{Fredholmness and essential spectrum}\label{sec6.2}

It is easy to see, using equivalence after extension and the theory of Wiener-Hopf factorization (WH factorization) that $ T_{z - \lambda}^M$ is Fredholm for all $ \lambda \in \mathbb{C} \setminus \mathbb{T}$. Indeed, this is a particular case of the following more general result. Here we denote by $\GR$ the space of all rational functions without poles on $\mathbb{T}$.
\begin{thm}\label{thm6.12}
Let $R$ be a rational function without zeros or poles on $\mathbb{T}$, i.e., $R \in \mathcal{G}\GR$. Then $A_R^M$ is Fredholm.
\end{thm}

\beginpf
It is enough to prove that $T_{G_R}$, with 
\[
G_R = \begin{pmatrix}
\overline\theta & 0 & 0 & 0 \\
0 & \overline\theta & 0 & 0 \\
R & R \overline{A_+} \overline{\theta} & \theta & 0 \\
R A_- \overline{\theta} & R & 0 & \theta 
\end{pmatrix}
\]
is Fredholm. This follows from the fact that $G_R$ admits a meromorphic factorization (\cite {CLS, CM}, see also \cite{CLS1}, Theorem 3.3) of the form $G_R = M_- M_+^{-1}$ with $M_-^{\pm1} \in (\overline{H^{\infty}}+ \GR)^{4\times 4}$
and $M_+^{\pm1} \in  ({H^{\infty}}+ \GR)^{4\times 4}$ given by:
\begin{eqnarray*}
M_+ &=& \begin{pmatrix}
1 & 0 & \theta & 0 \\
0 & 1 & 0 & \theta \\
0 & 0 & -R & 0 \\
0 & 0 & 0 & -R
\end{pmatrix}, \nonumber
\\
M_- &=& \begin{pmatrix}
\overline\theta & 0 & 1 & 0 \\
0 & \overline\theta & 0 &1 \\
R & R\overline{A_+}\ol\theta & 0 & \ol{A_+}R \\
RA_-\ol\theta & R & A_-R & 0
\end{pmatrix}.
\end{eqnarray*}
as can be easily verified.
\endpf

\begin{cor}\label{cor6.13}
$\sigma_e (T_z^M) \subset \mathbb{T}$.
\end{cor}

Now, to study the Fredholmness of $A_{z-\lambda}^M$ for $ \lambda \in \mathbb{T}$ we use another factorization of $\G_\lambda$. Indeed, more generally, for any $R\in \GR$ we can factorise
\[
G_R = \underbrace{\begin{pmatrix}
I_{2 \times 2} & 0_{2 \times 2} \\
R \begin{pmatrix}
1 & \overline{A_+} \\
A_- & 1
\end{pmatrix} & I_{2 \times 2}
\end{pmatrix}}_{H_R^-}
\underbrace{\begin{pmatrix}
\overline{ \theta} I_{2 \times 2} & 0_{2 \times 2} \\
R ( 1 - \overline{\theta}) I_{2 \times 2} & \theta I_{2 \times 2}
\end{pmatrix}}_{\widetilde{G}_R}
\]
where both factors depend on the symbol $ R \in \mathcal{R}$ (or the point $\lambda$ if $R = z - \lambda$), but the roles of $A_-, A_+$ ( i.e., $ \phi$ and $\psi$) and $ \theta$ are separated. 

Since $H_R^- \in \mathcal{G} ( \overline{H^{\infty}} + \mathcal{R})^{4 \times 4}$, it follows that ${G}_R$ admits a WH-factorization in $L^2$ \cite{CDR} (also known as generalised factorisation \cite{MP} or $\Phi$-factorisation \cite {LS}), whose existence is equivalent to $T_{G_R}$ being Fredholm \cite{CDR}, if and only if $\widetilde{G}_R$ admits such a factorization (Theorem 3.10 in \cite{MP}). Since $\widetilde{G}_R$ does not depend on $A_+$ or $A_-$, and for $A_+ = 0$ or $A_- = 0$ the operator $W$ defined in Theorem \ref{thm:eqblock} is triangular, taking into account this theorem we conclude the following:
\begin{thm}\label{thm6.14}
$T_R^M$ is Fredholm if and only if $A_R^{\theta}$ is Fredholm.
\end{thm}
The Fredholmness of truncated Toeplitz operators with rational symbols $ R \in \mathcal{R}$ was studied in Section 5 of \cite{CP16} in the equivalent setting of the real line. For the case $R= z - \lambda$, taking into account also Corollary \ref{cor6.13}, we have:
\begin{cor}\label{cor6.15}
The essential spectrum $\sigma_e ( T_z^M)$ is contained in $ \mathbb{T}$ and does not depend of $A_+$ or $A_-$, but only on the point $ \lambda$ and the inner function $\theta$. We have 
\[
\sigma_e ( T_z^M) = \sigma_e ( A_z^{\theta}) = \sigma(\theta) = \{ \lambda \in \mathbb{T} : \liminf_{z \to \lambda, z \in \mathbb{D}} | \theta (z) | = 0 \}.
\]

\end{cor}
\begin{rem}
Since, by Theorems \ref{4.7} and \ref{thm6.8}, the dimensions of $ \ker T_{z - \lambda}^M$ and $ \ker (T_{z - \lambda}^M)^* $ are equal and finite, we see that if $\lambda \in \sigma ( \theta)$ the range of $T_{z - \lambda}^M$ is not closed.
\end{rem}

\subsection{Invertibility, spectrum and resolvent operators}
From Corollary \ref{cor6.15} and from the description of $ \ker T_{z - \lambda}^M$ obtained in Section \ref{sec6.1}, taking moreover into account Corollary \ref{cor:n3.8}, we easily get the spectrum of the dual-band  shift.

\begin{thm}\label{thm6.17}
\[
(i) \hskip 0.1cm \lambda \in  \sigma(T^M_z) \iff (\lambda \in \mathbb{D}, \Dk =0) \vee (\lambda \in \mathbb{D}^- , \widetilde{\Dk} = 0 ) \vee ( \lambda \in \mathbb{T} \cap \sigma ( \theta)) 
\]
\[
\vee ( \lambda \in \mathbb{T}_{\text{ADC}}\setminus \sigma ( \theta ), \Dk = 0 ) ;
\]
(ii)
\begin{align*}
\hskip 0.1cm \sigma_p (T_z^M)& = \{ \lambda \in \mathbb{D} , \Dk =0 \} \cup \{ \lambda \in \mathbb{D}^- , \widetilde{\Dk}=0 \} \cup \{ \lambda \in \mathbb{T}_{\text{ADC}} , \Dk = 0 \}, \\
\sigma_e (T_z^M)& = \mathbb{T} \setminus \mathbb{T}_{\text{ADC}} \cup \{ \lambda \in \mathbb{T}_{\text{ADC}} \cap \sigma ( \theta), \Dk \ne 0 \}\\
\sigma_r ( T_z^M)& = \emptyset
\end{align*}
where $\Dk $ and $ \widetilde{\Dk}$ were defined in Propositions \ref{prop6.1} and \ref{prop6.2}.
\end{thm}

From Theorem \ref{thm6.17} we see in particular that unlike the essential spectrum, $\sigma ( T_z^M )$ is in general clearly different from $\sigma  ( A_z^{\theta})$. 

\vspace{2mm}

For $ \lambda \notin \sigma ( T_z^M)$, i.e., $ \lambda \in \mathbb{D},$ $\Dk \ne 0$, or $ \lambda \in \mathbb{D}^- $, $\widetilde{\Dk} \ne 0$, or $ \lambda \in \mathbb{T}_{\text{ADC}} \setminus \sigma ( \theta)$, $\Dk \ne 0$, we can explicitly define the resolvent operator by using Corollary \ref{4.6} and a bounded canonical Wiener-Hopf factorization of $\G_\lambda$, of the form $\G_{\lambda} = \G_{\lambda -} \G_{\lambda +}$, since we have, in that case, $(T_{\G_\lambda})^{-1} = (\G_{\lambda +})^{-1} P^+ (\G_{\lambda -})^{-1} P^+_{|(H^2)^4}$ \cite{CDR}. That canonical factorization will be given below for $ \lambda \in \mathbb{D} \cup \mathbb{T}_{\text{ADC}} \setminus \sigma ( \theta)$ with $ \Dk \ne 0$; for $ \lambda \in \mathbb{D}^-$, $\widetilde{\Dk} \ne 0$, the canonical factorization can be obtained analogously. Those factorizations, which were obtained by solving a Riemann-Hilbert problem of the form $\G_\lambda \phi_+ = \phi_-$ with $ \phi_+ \in (H^2)^4$ and $ \phi_- \in (\overline{H^2})^4$ for each of the column factors as in \cite{CCMN}, can be checked directly by multiplication of the matricial factors.

Note that if $ \lambda \in \mathbb{T} \setminus \sigma ( \theta )$ then $ \theta$ has an analytic continuation to a neighbourhood of $ \lambda$ and there is clearly an ADC for $ \theta$ there with $\frac{\theta - \theta( \lambda)}{z - \lambda} \in H^{\infty}$. If $ \lambda \in \mathbb{T}_{\text{ADC}} \cap \sigma ( \theta ) $ then $\frac{\theta - \theta( \lambda)}{z - \lambda}  \in H^2 \setminus H^{\infty}$ (page 505 of \cite{sarason07}); in this case $ \G_\lambda$ has an $L^2$-factorization \cite{LS}, but it is not bounded nor a WH factorization, although the factors are given by the same expressions as in the theorem below.

\begin{thm}\label{thm6.18}
(i) \hskip 0.1cm If $ \lambda \in \mathbb{D}$ or $ \lambda \in \mathbb{T}_{\text{ADC}} \setminus \sigma ( \theta )$ and $\Dk \ne 0$ then $\G_\lambda$ admits a bounded canonical factorization of the form $\G_{\lambda} = \G_{\lambda -} \G_{\lambda +}$ where \\
(i) if $\Delta = 1 - \overline{A_+(0)}\overline{\overline{A_-}(0)} \overline{\theta(0)}^2 \neq 0,$
\newline
\[
\G_{\lambda + }^{-1}=
\begin{pmatrix}
\theta+\frac{ \ol{A_+(0)} \,\ol{\ol{A_-}(0)}\,\overline{\theta(0)}}{\Delta} & \frac{\theta-\theta(\lambda)}{z-\lambda} & \frac{-\ol{A_+(0)}}{\Delta} & 0 \\ \\
-\frac{\ol{\ol{A_-}(0)}}{\Delta} & 0 & \theta+\frac{ \ol{A_+(0)}\, \ol{\ol{A_-}(0)}\,\overline{\theta(0)}}{\Delta} & \frac{\theta-\theta(\lambda)}{z-\lambda} \\ \\
-(z-\lambda) & -1 & 0 & 0 \\ \\
0 & 0 & -(z-\lambda) & -1
\end{pmatrix},
\]


and 

\[
\G_{\lambda - }=
\begin{pmatrix}
1+ \frac{ \ol{A_+(0)}\,\ol{\ol{A_-}(0)}\,\ol{\theta(0)} }{\Delta} \ol\theta & 
\frac{1-\theta(\lambda)\ol\theta}{z-\lambda}
& - \frac{\ol{A_+(0)}}{\Delta}\ol\theta & 0 \\ \\
- \frac{\ol{\ol{A_-}(0)}}{\Delta}\ol\theta & 0 & 1+ \frac{ \ol{A_+(0)}\,\ol{\ol{A_-}(0)}\,\ol{\theta(0)} }{\Delta} \ol\theta & \frac{1-\theta(\lambda)\ol\theta}{z-\lambda} \\ \\
g_{31}^- & -\theta(\lambda) & g_{33}^- & \ol{A_+}(1-\theta(\lambda)\ol\theta) \\ \\
g_{41}^- & A_-(1-\theta(\lambda)\ol\theta) & g_{43}^- & -\theta(\lambda)
\end{pmatrix},
\]

with
\begin{eqnarray*}
g^-_{31} &=& -\frac{\abar}{\Delta}(z-\lambda)(\ol{A_+}\ol\theta-\ol{A_+(0)}\ol{\theta(0)}),\\
g^-_{41} &=& \frac{z-\lambda}{\Delta} (A_- - \abar + \ol{A_+(0)}\,\abar \ol{\theta(0)}\,A_-\, (\ol\theta-\ol{\theta(0)}))\\\
g_{33}^- &=&  \frac{z-\lambda}{\Delta} (\ol{A_+}-\ol{A_+(0)} + \ol{A_+(0)}\, \abar\, \ol{\theta(0)}\, \ol{A_+}\,(\ol\theta-\ol{\theta(0)})),\\
g_{43}^- &=& -\frac{\ol{A_+(0)}}{\Delta}(z-\lambda) (A_-\ol\theta - \abar \ol{\theta(0)}).
\end{eqnarray*}
Note that in this case $\det \G_{\lambda + }^{-1} = \det \G_{\lambda - } =- \Dk / \Delta \in \mathbb{C} \setminus \{ 0 \} $.
\\

\noindent (ii) If $\Delta=0$, in which case $A_+(0)$, $\ol{A_-}(0)$, $\theta(0) \ne 0$ and
${A_+(0)}{\theta(0)} = \dfrac{1}{\overline{A_-}(0) \,\theta(0)}$, we have

\[
\G_{\lambda + }^{-1}=
\begin{pmatrix}
- \ol{A_+(0)}(1-\theta\ol{\theta(0)}) & \frac{\theta-\theta(\lambda)}{z-\lambda} & - \ol{A_+(0)}\ol{\theta(0)} & 0 \\ \\
\theta & 0 & 1 & \frac{\theta-\theta(\lambda)}{z-\lambda} \\ \\
-\ol{\theta(0)}\ol{A_+(0)}(z-\lambda) & -1 & 0 & 0 \\ \\
-(z-\lambda) & 0 & 0 & -1
\end{pmatrix},
\]
and

{\footnotesize
\[
\G_{\lambda -}=
\begin{pmatrix}
-\ol{A_+(0)}(\ol\theta-\ol{\theta(0)}) & \frac{1-\theta(\lambda)\ol\theta}{z-\lambda} & -\ol{A_+(0)}\,\overline{\theta(0)}\, \ol{\theta} & 0 \\ \\
1 & 0 & \ol\theta & \frac{1-\theta(\lambda)\ol\theta}{z-\lambda} \\ \\
(\ol{A_+}-\ol{A_+(0)})(z-\lambda) & -\theta(\lambda) & (\ol{A_+}\ol\theta - \ol{A_+(0)}\ol{\theta(0)})(z-\lambda)
& \ol{A_+}(1-\theta(\lambda)\ol\theta) \\ \\
-\ol{A_+(0)}A_-(\ol\theta-\ol{\theta(0)})(z-\lambda) & A_- (1-\theta(\lambda)\ol\theta)
& 1-\ol{A_+(0)}\ol{\theta(0)}A_-\ol\theta(z-\lambda) & -\theta(\lambda)
\end{pmatrix}
\]
}
with 
\[
\det (\G_{\lambda +})^{-1}=\det\G_{\lambda - }=
\ol{A_+(0)}(1-2\ol{\theta(0)}\theta(\lambda))=-\frac{\Dk}{\overline{\overline{A_-}(0)}} \in \CC\setminus\{0\}.
\]
\end{thm}

\section{$L^2$-factorization and angular derivatives}

In section \ref{sec6.2} the Fredholmness of $T_{z - \lambda}^M$, or $T_{\G_{\lambda}}$, was studied from a factorization of $\G_\lambda$, which we repeat here for convenience:
\begin{equation}
    \G_\lambda = H_- \widetilde{\G}_{\lambda}
\end{equation}
with
\begin{equation}
    H_-^{\pm} \in ( \overline{H^{\infty}} + \mathcal{R})^{4 \times 4} \text{ and } \widetilde{\G}_{\lambda} = \begin{pmatrix}
\overline{ \theta} I_{2 \times 2} & 0_{2 \times 2} \\
(z - \lambda) ( 1 - \overline{\theta}) I_{2 \times 2} & \theta I_{2 \times 2}
\end{pmatrix},
\end{equation}
using the fact that $\G_\lambda$ admits a WH factorization if and only if $\widetilde{\G}_{\lambda}$ admits a factorization of the same type. Indeed it is well known that the existence of a WH factorization for $ G \in (L^{\infty})^{n \times n }$ is equivalent to the Fredholmness of the Toeplitz operator $T_G$ (see, e.g. \cite[Thm.~1.1]{CDR} or \cite{LS,MP}).

It may happen that $\G_{\lambda}$ admits an $L^2$-factorization which does not satisfy the condition of boundedness for the densely defined operator $\G_{\lambda +}^{-1} P^+ \G_{\lambda_-}^{-1}: \mathcal{R}^4 \to (L^2)^4$ (see \cite{CDR, LS})
Although in that case one cannot conclude that $T_{\G_{\lambda}}$ is Fredholm, it is still possible to use it to characterize several important properties of the operator, such that as injectivity (\cite{LS}). This leads us to the question of existence of such a $L^2$-factorization for $\widetilde{\G}_{\lambda}$. Somewhat surprisingly, we obtain a necessary and sufficient condition for existence of an ADC of $ \theta$ at $ \lambda \in \mathbb{T}$ in terms of an $L^2$-factorization, a relation which appears here for the first time.

\begin{thm}\label{thm7.1}
Let $ \lambda \in \mathbb{T}$. The inner function $\theta$ has an ADC at $\lambda$ if and only if $\widetilde{\G}_{\lambda}$ has an $L^2$ factorization.
\end{thm}
\begin{proof}
(i) Assume that $\widetilde{\G}_{\lambda} = \widetilde{\G}_{\lambda - } D \widetilde{\G}_{\lambda +}$ with $D = \diag \{ z^{k_j} \}$ with $k_j \in \mathbb{Z}$, $j = 1,2,3,4$. Since $\det \widetilde{\G}_{\lambda}= 1$, we must have $\sum_{j=1}^4 k_j =0$. If, for some $j$, we have $k_j = -n <0$, then there exists a non-zero solution to 
\begin{equation}\label{7.3}
    \begin{pmatrix}
\overline{ \theta} I_{2 \times 2} & 0_{2 \times 2} \\
(z - \lambda) ( 1 - \overline{\theta}) I_{2 \times 2} & \theta I_{2 \times 2}
\end{pmatrix} \begin{pmatrix}
\psi_{1+} \\
\psi_{2+}
\end{pmatrix} = \overline{z}^n  \begin{pmatrix}
\psi_{1-} \\
\psi_{2-}
\end{pmatrix} 
\end{equation}
with $\psi_{1+}, \psi_{2+} \in (H^2)^2$ and $
\psi_{1-}, \psi_{2-} \in (\overline{H^2})^2$, given by the $j$'th column of the factors. So we have, with $ \overline{\theta} \psi_{1+} = \overline{z}^n \psi_{1-}$,
\begin{align}
    &(z - \lambda) ( 1 - \overline{\theta}) \psi_{1+} + \theta \psi_{2+} = \overline{z}^n \psi_{2-} \nonumber \\
    \implies& (z - \lambda) \psi_{1+} - ( z - \lambda) \overline{z}^n \psi_{1-} + \theta \psi_{2+} = \overline{z}^n \psi_{2-} \nonumber \\
    \implies& (z - \lambda) \psi_{1+} + \theta \psi_{2+} = \overline{z}^n\psi_{2-} + (z - \lambda) \overline{z}^n \psi_{1-} = C \in \mathbb{C}^2 \label{7.4}.
\end{align}
If $C=0$ (which necessarily happens if $n >1$) then $\psi_{1+} = \psi_{2+} = 0$ and there is no nontrivial solution to \eqref{7.3}. So we must have $ n \leqslant 1$ and $C \ne 0$. For $ n \geqslant 1$, from \eqref{7.4},
\[
\overline{z}^n ( z - \lambda) \psi_{1-} - \overline{\theta}C= - \psi_{2+} = A ( z - \lambda) + B, \text{ with } A, B \in \mathbb{C}^2.
\]
Now we must also have $B \ne 0$ because, otherwise we would have $ \psi_{1-} - A = \overline{\theta} C \frac{z^n}{z - \lambda} \in \overline{H^2}$ with $ C \ne 0$, which is impossible because $ \lambda \in \mathbb{T}$. Going back to \eqref{7.4}, then,
\begin{align}
    &( z- \lambda) \psi_{1+} + \theta ( - A( z - \lambda) - B ) = C \nonumber \\
    \iff& (z - \lambda) (\psi_{1+} - \theta A ) = C + \theta B ; \nonumber 
\end{align}
by Lemma \ref{lem6.3}, it follows that $ \theta( \lambda)$ exists and, since $\psi_{1+} -  \theta A \in (H^2)^2$ we get that $ \frac{\theta - \theta( \lambda)}{z - \lambda} \in H^2$ and $ \theta $ has an ADC at $ \lambda$ by Lemma \ref{lem6.4}.
\\

\noindent (ii) Conversely, suppose that $ \theta $ has an ADC at $ \lambda$. Then $\widetilde{\G}_{\lambda}$ admits the following $L^2$ factorization (obtained as before by solving a Riemann-Hilbert problem):
\\
- if $ \theta( \lambda) \ne \frac{-1}{1 - \overline{\theta(0)}}$, $\widetilde{\G}_{\lambda} = \widetilde{\G}_{\lambda-} \widetilde{\G}_{\lambda +}$ with
{
\[
\widetilde{\G}_{\lambda-} =
\begin{pmatrix}
\frac{\overline{\theta}}{1 - \overline{\theta(0)}}+1 & \frac{1-\theta(\lambda)\ol\theta}{z-\lambda} & 0  & 0 \\ \\
0 & 0 & \frac{\overline{\theta}}{1 - \overline{\theta(0)}}+1 & \frac{1-\theta(\lambda)\ol\theta}{z-\lambda} \\ \\
- \frac{\overline{\theta} - \overline{\theta(0)}}{1 - \overline{\theta(0)}}( z - \lambda) & - 1 + \theta( \lambda) (\overline{\theta}-1) & 0
& 0 \\ \\
0 & 0 & - \frac{\overline{\theta} - \overline{\theta(0)}}{1 - \overline{\theta(0)}}( z - \lambda) & - 1 + \theta( \lambda) (\overline{\theta}-1)
\end{pmatrix}
\]
}

{
\[ 
\widetilde{\G}_{\lambda +}^{-1} =
\begin{pmatrix}
\frac{1}{1 - \overline{\theta(0)}} +  \theta & \frac{\theta - \theta( \lambda)}{z - \lambda} & 0  & 0 \\
0 & 0 & \frac{1}{1 - \overline{\theta(0)}} +  \theta & \frac{\theta - \theta( \lambda)}{z - \lambda} \\ 
- ( z - \lambda) & -1 & 0
& 0  \\
0 & 0 & - ( z - \lambda) & -1
\end{pmatrix},
\]
}
\\

\noindent -if $ \theta( \lambda) = \frac{-1}{1 - \overline{\theta(0)}}$, $\widetilde{\G}_{\lambda} = \widetilde{\G}_{\lambda-} \diag{( \overline{z}, \overline{z},z,z)} \widetilde{\G}_{\lambda +}$ with
{
\[ 
\widetilde{\G}_{\lambda - } =
\begin{pmatrix}
z \frac{1-\theta(\lambda)\ol\theta}{z-\lambda}  & 0  & \overline{z} & 0 \\
0 & z \frac{1-\theta(\lambda)\ol\theta}{z-\lambda} & 0 & \overline{z} \\ 
z ( \theta ( \lambda ) \overline{\theta} -1 - \theta( \lambda)) & 0 & -\frac{z - \lambda}{z}
& 0  \\
0 & z ( \theta ( \lambda ) \overline{\theta} -1 - \theta( \lambda)) & 0 & -\frac{z - \lambda}{z}
\end{pmatrix},
\]
}
{
\[ 
\widetilde{\G}_{\lambda + }^{-1} =
\begin{pmatrix}
\frac{\theta - \theta(\lambda)}{z - \lambda}  & 0  & \theta & 0 \\
0 & \frac{\theta - \theta(\lambda)}{z - \lambda} & 0 & \theta \\ 
-1 & 0 & -(z - \lambda)
& 0  \\
0 & -1 & 0 & -(z - \lambda)
\end{pmatrix}.
\]
}
\end{proof}

\end{document}